\documentclass[10pt]{amsart}
\usepackage{amsthm}
\usepackage[english]{babel}
\input xy
\xyoption{all}
\usepackage{amssymb,amsmath,color}

\newcommand{\abs}[1]{\left\vert#1\right\vert}
\newcommand{\R}{\mathbb{R}}

\newcommand{\hook}{\lrcorner\,}

\newcommand{\SU}{\mathrm{SU}}

\newcommand{\Gtwo}{\mathrm{G}_2}

\newcommand{\GL}{\mathrm{GL}}

\newcommand{\dfn}[1]{\emph{#1}}

\DeclareMathOperator{\Span}{Span}

\theoremstyle{plain}

\theoremstyle{definition}

\theoremstyle{remark}
\newtheorem*{remark}{Remark}

\usepackage{moreverb}
\usepackage{graphicx}
\newcommand{\cpp}[1]{{\tt #1}}
\newcommand{\CC}{C\nolinebreak\hspace{-.05em}\raisebox{.4ex}{\tiny\bf +}\nolinebreak\hspace{-.10em}\raisebox{.4ex}{\tiny\bf +}}

\begin{document}
\sloppy

\title{Symbolic computations in differential geometry}
\author{Diego Conti}
\address{Dipartimento di Matematica e Applicazioni, Universit\`a di Milano Bicocca, via Cozzi 53, 20125 Milano, Italy.}
\email{diego.conti@unimib.it}


\begin{abstract}
We introduce the {\CC} library Wedge, based on GiNaC, for symbolic computations in differential geometry. We show how Wedge makes it possible to use the language {\CC} to perform such computations, and illustrate some advantages of this approach with explicit examples. In particular, we describe a short program to determine whether a given linear exterior differential system  is involutive.
\end{abstract}
\maketitle

\centerline{\small\textbf{Keywords}: Curvature, differential forms, computer algebra, exterior differential systems.}

\section*{Introduction}
There are many computationally intensive problems in differential and Riemannian geometry that are best solved by the use of a computer. Due to the nature of these problems, any system meant to perform this type of calculations must be able to manipulate algebraic and differential expressions. For this reason, such systems are generally implemented as extensions, or packages, for a general purpose computer algebra system: the latter takes care of handling expressions, and the extension introduces the differential-geometry specific features --- differential forms, tensors, connections and so on. Examples are given by the packages difforms and GRTensor (see \cite{GRTensor}) for Maple, or the Ricci package for Mathematica.

A remarkable consequence of this approach is that one is essentially limited, when implementing one's own algorithms, to using the programming language built in the computer algebra system. Some drawbacks and limitations of these programming languages are described in \cite{GiNaC}, where an alternative is also introduced, namely the {\CC} library GiNaC. As suggested by the name (an acronym for GiNaC Is Not A CAS), GiNaC differs from the above mentioned computer algebra systems in that it is based on a general purpose, well-established programming language such as \CC, rather than introducing a new one. In particular, development, debugging and documentation of a program based on GiNaC can take advantage of the many tools commonly available to a {\CC} programmer. This makes GiNaC a natural choice when implementing new, complex algorithms, like the one introduced in \cite{Conti:InvariantForms}, which provided the original motivation for the present work.

\smallskip
In this paper we introduce  Wedge, an extension to GiNaC that can be used to write a {\CC} program that performs computations in differential and Riemannian geometry. Wedge is able to perform algebraic or differential computations with differential forms and spinors, as well as curvature computations in an adapted frame, and contains some support for vector spaces, represented in terms of a basis. Bases of vector fields, or frames, play a central r\^ole in Wedge: in particular, the tangent space of a manifold is represented by a frame, and a Riemannian metric is represented by a (possibly different) orthonormal frame. Notice that the above-mentioned package  GRTensor also supports working with adapted frames, but our approach differs in that the frame need not be defined in terms of coordinates; this can be useful when working on a Lie group, where the geometry is defined naturally by the structure constants, or on the generic manifold with a fixed geometric structure (see Section~\ref{sec:inheritance}).

Another unique feature of Wedge among packages for differential and Riemannian geometry, beside the choice of  {\CC}, is the fact that it is completely based on free, open-source software. Like GiNaC, Wedge is licensed under the GNU general public license; its source code is available at  {\tt http://libwedge.sourceforge.net}.

\smallskip
This paper is written without assuming the reader is familiar with {\CC}, and its purpose is twofold: to introduce the main functionality of Wedge, and to illustrate with examples certain features of {\CC} which can prove very helpful in the practice of writing a program to perform some specific computation.

In the first section we introduce some basic functionality, concerning differential forms and connections; at the same time, we illustrate classes and inheritance.

In the second section we introduce spinors and ``generic'' manifolds; at the same time, we illustrate object-oriented programming.

In the third section we explain briefly how GiNaC handles expressions, then introduce bases and frames, namely the linear algebra features of Wedge.

In the final section we give an application from Cartan-K\"ahler theory, with a short program that reproduces the computations of \cite{Bryant} that prove the local existence of metrics with holonomy $\Gtwo$.

\section{Working with differential forms and connections}
\label{sec:classes} A fundamental feature of {\CC} is the possibility of introducing user-defined
types, i.e. classes (or structs). The definition of a class specifies not only the type of data
contained in a variable which has that class as its type, but also some operations that can be
performed on this data. This is generally better than having global functions which either take a
long list of arguments or use global variables, and this can be seen, for instance, in situations where
more sets of data appear in the same program. In this section we shall illustrate this point with
an example, considering a problem concerning multiple connections on a fixed manifold. In the
course of the section, we shall introduce some of the essential functionality of Wedge.

The basic geometric entity in Wedge is the class \cpp{Manifold}. A variable of type \cpp{Manifold}
represents a manifold in the mathematical sense, which is assumed to be parallelizable, and
represented by a global basis of one-forms $e^1,\dotsc,e^n$, with dual basis of vector fields $e_1,\dotsc, e_n$. This assumption is tailored on Lie
groups, but one can also think of a \cpp{Manifold} object as representing a coordinate patch;
also, every manifold can be written as a quotient $M/G$ of a parallelizable manifold, so one can always reduce to the parallelizable case. For instance, one can compute the curvature of a Riemannian metric on $M/G$ applying the O'Neill formula (see \cite{ONeill}).

As an example, we shall
consider the nilpotent Lie group $X$, characterized by the existence of a global basis of one-forms
$e^1,\dotsc, e^4$ such that
\begin{equation}
\label{eqn:nilmanifold}de^1=0,\quad de^2=0,\quad de^3=e^{12},\quad de^4=e^{13}.
\end{equation}
In the above equation, $e^{23}$ stands for the wedge product \mbox{$e^2\wedge e^3$}, and so on;
this shorthand notation is also used by Wedge, and so it will appear throughout the paper in both formulae and program output.
Equations~\eqref{eqn:nilmanifold} translate to the following code:
\smallskip\begin{verbatimtab}
struct X : public ConcreteManifold , public Has_dTable {
    X() : ConcreteManifold(4) {
        Declare_d(e(1),0);
        Declare_d(e(2),0);
        Declare_d(e(3),e(1)*e(2));
        Declare_d(e(4),e(1)*e(3));
    }
};
\end{verbatimtab}

\begin{figure}
\label{fig:ManifoldInheritanceDiagram}
\includegraphics[width=\textwidth,height=5.5cm]{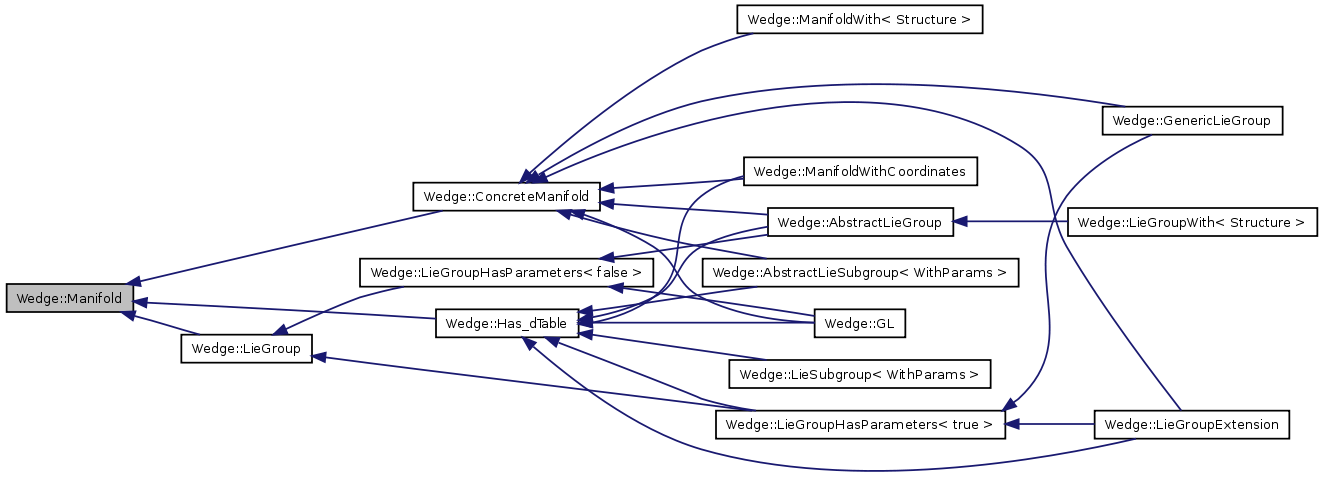}
\caption{Inheritance graph for classes inheriting from \cpp{Manifold}}
\end{figure}

The first line means that \cpp{X} \emph{inherits} from the two classes \cpp{ConcreteManifold} and
\cpp{Has\_dTable}. Roughly speaking, this means that it inherits the functionality implemented by
these two classes; inheritance defines a partial ordering relation, and is best represented by a
graph (see Figure~\ref{fig:ManifoldInheritanceDiagram}).
Specifically, inheriting from \cpp{ConcreteManifold} ensures that the forms $e^1,\dotsc,e^4$ are
defined (where the dimension $4$ appears as a parameter on the second line). The fact that \cpp{X}
inherits from \cpp{Has\_dTable} means that the operator $d$ on the manifold $X$ is known in terms
of its action on the basis $e^1,\dotsc,e^4$; the relations \eqref{eqn:nilmanifold} are given with
calls to the member function \cpp{Declare\_d}. These calls appear in the body of the
\emph{constructor} of \cpp{X}, which is invoked automatically when a variable of type \cpp{X} is
constructed. We can now instantiate a variable of type \cpp{X} and perform computations with it,
e.g.
\begin{verbatimtab}
    X M;
    cout<<M.d(M.e(4))<<endl;
\end{verbatimtab}
has the effect of printing $e^{13}$. Notice that both the forms $e^i$ and the operator $d$ are implemented as
\emph{members} of \cpp{X}; this means that every variable of type \cpp{X} has its own set of forms
$e^1,\dotsc,e^n$ and $d$ operator, and so the code must tell the compiler which \cpp{X} should be
used. This is achieved here by the ``\cpp{M.}'' appearing in the function calls; however, inheritance
provides an equivalent, more appealing alternative, as in the following:
\begin{verbatimtab}
struct SomeCalculations : X {
    SomeCalculations() {cout<<d(e(4))<<endl;}
} M;
\end{verbatimtab}

\smallskip
The class \cpp{Manifold} also implements the Lie bracket and Lie derivative as member functions,
whilst the exterior and interior product of forms are implemented by global functions. Covariant
differentiation and curvature computations are accounted for by the class \cpp{Connection}, which we
now introduce with an example.

Every almost-complex manifold  $(M,J)$ admits an almost complex connection $\tilde\nabla$ such
that its torsion form $\Theta$ and the  Nihenjuis tensor are related by
\[16\Theta(X,Y)=N(X,Y)\;;\]
by \cite{KobayashiNomizu} $\tilde\nabla$ can be obtained from an arbitrary torsion-free connection
$\nabla$ by
\begin{equation}
 \label{eqn:tildenabla}
\tilde\nabla_XY=\nabla_X Y-Q(X,Y),\quad 4Q(X,Y)=(\nabla_{JY}J)X+J((\nabla_Y J)X)+2J((\nabla_X J)Y).
\end{equation}
Suppose one wants to compute such a connection in the case of our nilpotent Lie group $X$, with
the almost-complex structure determined by
 \[J(e_1)=e_2,\quad J(e_3)=e_4.\]
We can do so with the following
code:
\smallskip\begin{verbatimtab}
struct AlmostComplex : public X {
    TorsionFreeConnection<true> h;
    ex J(ex Y) {return Hook(Y,e(1)*e(2)+e(3)*e(4));}
    ex A(ex X,ex Y) {
        return h.Nabla<VectorField>(X,J(Y))-J(h.Nabla<VectorField>(X,Y));
    }
    ex Q(ex X, ex Y) {
        return (A(J(Y),X)+J(A(Y,X))+2*J(A(X,Y)))/4;
    }
    AlmostComplex() : omega(this,e()) {
        Connection k(this,e());
        for (int i=1;i<=4;++i)
        for (int j=1;j<=4;++j)
            k.DeclareNabla<VectorField>(e(i),e(j),
                h.Nabla<VectorField>(e(i),e(j))-Q(e(i),e(j)));
        cout<<k<<k.Torsion();
    }
} M;
\end{verbatimtab}

This code defines and instantiates a struct named \cpp{AlmostComplex} inheriting from \cpp{X};
\cpp{AlmostComplex} has a data member  \cpp{h} of type \cpp{TorsionFreeConnection<true>}, which
represents a generic torsion-free connection $\nabla$ on $X$. The constructor uses $h$ to compute
another connection \cpp{k}, corresponding to $\tilde\nabla$ in \eqref{eqn:tildenabla}, then prints
the connection forms of $\tilde\nabla$ and its torsion form $\Theta$. The latter is a vector
\[\begin{pmatrix} 0,&0,&-\frac{1}{4} e^{32}-\frac{1}{4} e^{41},&\frac{1}{4} e^{42}-\frac{1}{4}
e^{31}\end{pmatrix}\]
consistently with the fact that the almost complex structure $J$ is not
integrable.

There are two connection objects appearing here, both represented internally by their connection
forms in terms of the frame $e_1,\dotsc,e_n$. On construction, they are initialized in terms of
symbols $\Gamma_{ijk}$, representing the functions
\[\Gamma_{ijk}=\langle \nabla_{e_i} e_j,e^k\rangle,\tilde\Gamma_{ijk}=\langle \tilde\nabla_{e_i} e_j,e^k\rangle,\]
where $\langle,\cdot \rangle$ represents the pairing on $TM\otimes T^*M$. Since the type of \cpp{h} is \cpp{TorsionFreeConnection<true>}, the constructor of $\cpp{h}$
ensures that the torsion is zero by solving the linear equations in the $\Gamma_{ijk}$ given by
\begin{equation}
\label{eqn:TorsionFree}
\sum_i e^i\wedge \nabla_{e_i}e^j-de^j=0.
\end{equation}
Of course this does not determine $\nabla$ uniquely, and some of the $\Gamma_{ijk}$ remain as
parameters. The conditions \eqref{eqn:tildenabla} on the connection $k$ are imposed directly by
calling the member function \cpp{DeclareNabla}. The tensor $Q$ is computed by the function
\cpp{Q}, which in turn calls \cpp{Nabla}; since the latter is a member function, the connection $h$ is
used to compute the covariant derivative. Thus, the connection forms of $k$ depends on the
$\Gamma_{ijk}$, although the torsion $\Theta$ is independent of the $\Gamma_{ijk}$, as expected
from the theory.

Another important point is that \cpp{DeclareNabla} does not assign values directly to the
connection forms. Rather, it uses GiNaC's function \cpp{lsolve} to solve a system of equations,
linear in the $\Gamma_{ijk}$, and substitutes the solution into the connection forms. This
procedure is more flexible than giving the list of the $\Gamma_{ijk}$ for two reasons: the
arguments passed to \cpp{DeclareNabla} need not be elements of the frame (they can even be forms
of degree higher than one, or, for Riemannian connections only, spinors), and part of the
connection forms may be left unspecified. However, said procedure only works because the
$\Gamma_{ijk}$ are stored in the \cpp{Connection} object, which is how \cpp{DeclareNabla} knows it
should leave the $\Gamma_{ijk}$ as parameters, and solve with respect to the $\tilde\Gamma_{ijk}$.

Summing up, we have illustrated how the fact that a connection is represented as an instance of a class makes working with two different connections on a fixed manifold as natural as defining
two variables of the same type.

\section{Torsion free connections and generic manifolds}
\label{sec:inheritance} Class \cpp{X} from Section~\ref{sec:classes} represents a fixed manifold,
where the action of the $d$ operator can be recovered from its action on a basis of one-forms; we
have seen that Wedge enables one to perform torsion computations, or impose torsion conditions.
One can also go the other way, and compute the action of $d$ in terms of the covariant derivative,
with respect to a torsion-free connection. To illustrate this, suppose one has a four-dimensional
Riemannian manifold $X$, with a global orthonormal frame $e^1,\dotsc, e^4$. The parallelism
induces a spin structure; in general, on a  parallelizable $n$-dimensional manifold a spinor can
be viewed as a map $X\to\Sigma$, where $\Sigma$ is the $n$-dimensional spinor representation. Let  $m=[n/2]$ be the integer part of $n/2$;
$\Sigma$ is a complex vector space of dimension $2^m$, represented in Wedge in terms of the basis
$u_0,\dotsc,u_{2^{m}-1}$, where $u_k$ corresponds to
\[u((-1)^{\epsilon_{m-1}},\dotsc,(-1)^{\epsilon_{0}}),\quad k=\sum_{r=0}^{m-1} (1-\epsilon_r)2^r\] in the
notation of \cite{Seminarbericht} (see also \cite{ContiFino} for more explicit formulae). In our
case, we can declare that the constant spinor $\psi\colon X\to\Sigma$, $\psi\equiv u_0$ is
parallel, i.e.
\[\nabla_{e_i} \psi=0\,\quad i=1,\dots,4\]
 with the following code:
\smallskip\begin{verbatimtab}
struct X : public ManifoldWith<RiemannianStructure> {
    X() : ManifoldWith<RiemannianStructure>(4) {
        for (int i=1;i<=4;++i)
	        DeclareNabla<Spinor>(e(1),u(0),0);
    }
};
\end{verbatimtab}

Hence, $X$ is a generic parallelizable Riemannian $4$-manifold with holonomy contained in
$\SU(2)$; this condition does not determine the connection form uniquely, but it does impose
certain conditions on them. Internally, the template class \cpp{ManifoldWith<RiemannianStructure>}
contains a member of type \cpp{LeviCivitaConnection<false>}, which behaves in a similar way to its
counterpart \cpp{TorsionFreeConnection<true>} from Section~\ref{sec:classes}, with two differences: first, it represents the Levi-Civita connection, which is both torsion-free and Riemannian, and secondly it does
not use the operator $d$ of the manifold to impose the  condition $d^2=0$. Instead,
\cpp{ManifoldWith} uses the Levi-Civita connection to compute the action of $d$ via
\eqref{eqn:TorsionFree}.

\smallskip
\label{p:convention}
\textbf{Convention.} All further code fragments appearing in this section will be assumed to appear in
the body of the constructor of a class such as \cpp{SomeCalculations} of Section
\ref{sec:classes}.

In fact, whilst the new \cpp{X} does not derive from \cpp{has\_dTable}, the definition of
\cpp{SomeCalculations} of Section \ref{sec:classes} is still legitimate. Of course, the output
will depend on the $\Gamma_{ijk}$ which in part are left unspecified. In particular, it is not
possible to guarantee that $d^2=0$, as this leads to an underdetermined system of differential
equations, quadratic in the $\Gamma_{ijk}$ and linear in their derivatives. However, there are
many interesting relations not involving the derivatives of the $\Gamma_{ijk}$ that can be proved
using \cpp{ManifoldWith}. As a first example, we can verify that the existence of a parallel
spinor is equivalent to the existence of a local frame $e^1,\dotsc, e^4$ such that the forms
\begin{equation}
e^{12}+e^{34},\quad e^{13}+e^{42},\quad e^{14}+e^{23},
\label{eqn:TwoForms}
\end{equation}
are closed (see  \cite{Hitchin:SelfDuality}). In our case, since we have chosen $u_0$ as the
parallel spinor, the standard frame $e^1,\dotsc, e^4$ satisfies this condition, and  the code
\smallskip\begin{verbatimtab}
cout<<d(e(1)*e(2)+e(3)*e(4))<<endl;
cout<<d(e(1)*e(3)+e(4)*e(2))<<endl;
cout<<d(e(1)*e(4)+e(2)*e(3))<<endl;
\end{verbatimtab}
prints three times zero. One can obtain the opposite implication by invoking \cpp{Declare\_d} in the constructor of \cpp{X} to impose that the forms \eqref{eqn:TwoForms} are closed.

Compared to the examples of Section~\ref{sec:classes}, this code shows an important feature of {\CC} which is common to many programming languages, but not so common among computer algebra systems: one can define different functions with the same name, and the compiler is responsible for selecting the correct one based on the context. In this case the functions \cpp{d} and \cpp{Declare\_d} appear as members of different classes, i.e. \cpp{Has\_dTable} and \cpp{ManifoldWith}, and so the meaning of a call to \cpp{d} in the body of a member of \cpp{X} depends on whether \cpp{X} inherits from \cpp{Has\_dTable}, as was the case in Section~\ref{sec:classes}, or \cpp{ManifoldWith}, as in the example above. In fact, {\CC} allows even greater flexibility, as we now illustrate with a second example.

\smallskip
In order to give a slightly more complicated application, we observe that our manifold $X$ is in particular an almost-K\"ahler manifold with respect to any one of the closed two-forms \eqref{eqn:TwoForms}; we shall fix the first one. Thus, it makes sense to consider the bilagrangian splitting $TX=F\oplus G$, where $F=\langle e_1,e_3\rangle$ and $G=\langle e_2,e_4\rangle$;
this splitting determines a canonical connection $\omega$. By \cite{Vaisman:SymplecticCurvatureTensors}, the torsion of $\omega$ is zero if and only if the distributions $F$ and $G$ are integrable. We can actually prove this equivalence with Wedge; we shall illustrate the ``only if'' implication here.
\smallskip\begin{verbatimtab}
Connection omega(this,e(),"Gamma'");
for (int k=1;k<=4;++k) {
    omega.DeclareNabla<DifferentialForm>(e(k),e(1)*e(2)+e(3)*e(4),0);
    for (int i=1;i<=4;++i)
    for (int j=i%2+1;j<=4;j+=2)
        omega.DeclareZero(Hook(e(j),
            omega.Nabla<DifferentialForm>(e(k),e(i))));
}
for (int k=1;k<=4;++k)
for (int i=k%2+1;i<=4;i+=2)
for (int j=k%2+1;j<=4;j+=2)
    omega.DeclareZero(Hook(e(j),
        omega.Nabla<VectorField>(e(k),e(i))-LieBracket(e(k),e(i))));
exvector T=omega.Torsion();
DeclareZero(T.begin(),T.end());
cout<<LieBracket(e(1),e(3))<<","<<LieBracket(e(2),e(4))<<endl;
\end{verbatimtab}

This code defines a generic connection $\omega$, whose connection parameters are denoted by $\Gamma'_{ijk}$ to distinguish them from the parameters $\Gamma_{ijk}$ of the Levi-Civita connection, and imposes the two sets of conditions that determine the symplectic connection. First, the holonomy is reduced to $\GL(2,\R)$ by requiring that the symplectic form be parallel and the two Lagrangian distributions $F$ and $G$ preserved by the covariant derivative. Then one imposes
\[\nabla_{X_F} Y_G =[X_F,Y_G]_G, \quad \nabla_{X_G} [Y_F] =[X_G,Y_F]_F, \quad X,Y\in\Gamma(TX),\]
where the subscripts denote projection (reflected in the code by the use of the interior product function \cpp{Hook}). Notice that \cpp{LieBracket} is a  member function of \cpp{Manifold}, which \cpp{ManifoldWith} reimplements in terms of the Levi-Civita connection; thus, its result depends on the $\Gamma_{ijk}$. By the general theory, the conditions determine $\omega$ completely, i.e. $\omega$ does not depend on the $\Gamma'_{ijk}$ but only on the $\Gamma_{ijk}$. Since this dependence is linear, the condition that the torsion is zero gives equations in the $\Gamma_{ijk}$ that can be solved by the function \cpp{DeclareZero}, a member of \cpp{ManifoldWith}. Having imposed this torsion conditions, the program concludes that $[e_1,e_3]$ and $[e_2,e_4]$ equal
\[e_{1} \Gamma_{142}+ e_{3} \Gamma_{331},
\Gamma_{431} e_{4}+ \Gamma_{231} e_{2}\]
respectively, proving that the distributions $F$ and $G$ are involutive.

The code at work here  uses a typical object-oriented programming technique. Specifically, \cpp{Connection} needs to compute the action of $d$ in order to compute the torsion; to this end, every \cpp{Connection} object contains a pointer to the \cpp{Manifold} object it refers to. In this case, the manifold to which  $\omega$  refers is represented by a \cpp{ManifoldWith} object, which implements \cpp{d} using the Levi-Civita connection; since \cpp{d} is a \emph{virtual function}, the calls to \cpp{d} performed by \cpp{Torsion} execute the implementation of \cpp{ManifoldWith}. In the examples of Section~\ref{sec:classes}, those same calls execute the implementation of \cpp{Has\_dTable}. This ensures that the torsion of $\omega$ is computed correctly in terms of the $\Gamma_{ijk}$ here, whereas the code of Section~\ref{sec:classes} computes the torsion using the action of $d$ on the forms $e^1,\dotsc,e^4$, although the function \cpp{Torsion} is exactly the same.

\begin{remark} There is nothing essential about the assumption that $X$ has holonomy $\SU(2)$. In fact, one can easily modify the code to obtain the same result for a generic symplectic $4$-manifold $X$.
\end{remark}

\section{Forms and frames}
\label{sec:frames} In this section we introduce the main linear algebra functionality of Wedge;
before doing that, we need to explain how GiNaC handles expressions.

We have seen in Section~\ref{sec:classes} that manifolds are represented in Wedge by a global basis of one-forms $e^1,\dotsc, e^n$, which can be accessed via a member function of the class \cpp{Manifold}.
The {\CC} type of a form, say \cpp{X.e(1)}, is the class \cpp{ex} defined in the library GiNaC, which handles expressions such as \cpp{2*e(1)+e(2)}, the result of which is still an \cpp{ex}. The class \cpp{ex} acts as a proxy to the GiNaC class \cpp{basic}. In practice, this means that an \cpp{ex} contains the address of an object whose type is a class that inherits from \cpp{basic} (e.g. \cpp{add}, representing a sum, or \cpp{numeric} representing a number), and it is this object that performs the actual computations. The use of virtual functions ensures that the correct code is used, depending on the {\CC} type of the object the \cpp{ex} points to (for details, see \cite{GiNaC}).

Thus, while \cpp{X.e(1)} has type  \cpp{ex}, under the hood it ``refers'' to an object of type \cpp{DifferentialOneForm}, and it is this type that implements skew-commutativity, making use of a canonical ordering that exists among \cpp{basic} objects. More significantly,  objects of type \cpp{DifferentialOneForm} and linear combinations of them can be used interchangeably; nonetheless, they can be distinguished by their internal representation, since the former have type \cpp{DifferentialOneForm} and the latter have type \cpp{add} (e.g. \cpp{e(1)+e(2)}) or \cpp{mul} (e.g. \cpp{2*e(1)}). We refer to objects in the former category as \dfn{simple} elements. Mathematically, we are working with sparse representations of elements in the vector space generated by all the  variables of type \cpp{DifferentialOneForm} appearing in a program. All this works equally well with other types than  \cpp{DifferentialOneForm}.

In general, given a small set of vectors, it is possible to switch from sparse to dense representation by extracting a basis, and writing the other terms in components. This is accomplished by the class template \cpp{Basis}. As a simple example, suppose \cpp{X} defines the Iwasawa manifold, characterized by the existence of an invariant basis $e^1,\dotsc, e^6$ satisfying
\[de^1=0=de^2=de^3=de^4,\quad de^5=e^{13}+e^{42},\quad de^6=e^{14}+e^{23}.\]
We can compute a basis of the space of invariant exact three-forms with the following code (with
the usual convention of p.~\pageref{p:convention}):
\smallskip
\begin{verbatimtab}
Basis<DifferentialForm> b;
for (int i=1;i<=6;++i)
for (int j=i+1;j<=6;++j)
    b.push_back(d(e(i)*e(j)));
cout<<b;
\end{verbatimtab}
resulting in the output
\[e^{412},
- e^{312},
- e^{342},
e^{341},
e^{642}-e^{631}-e^{352}+e^{415}.\]
One can then recover the components of, say, $de^{45}$ by
\smallskip
\begin{verbatimtab}
cout<<b.Components(d(e(4)*e(5)));
\end{verbatimtab}
obtaining the vector
$\left(\begin{smallmatrix}0&0&0&-1&0\end{smallmatrix}\right)$.

\smallskip
Another use of \cpp{Basis} is related to the identification of each vector space with its dual. Indeed, the existence of simple elements enables us to establish a canonical pairing, which is bilinear and satisfies
\[\langle \alpha,\beta\rangle =\begin{cases} 1 & \alpha=\beta\\ 0 & \alpha\neq\beta\end{cases}\;,\quad \alpha,\beta\text{ simple.}\]
In particular, we can identify vector fields with one-forms, by interpreting the pairing $\langle,\rangle$ as the interior product. For instance, the standard frame $e^1,\dots, e^n$ associated to a \cpp{ConcreteManifold} always consists of simple elements, and so in this case the form $e^i$ and the vector field $e_i$ are represented by the same element \cpp{e(i)}. However, this fact only holds for $\langle ,\rangle$-orthonormal bases. In general, the dual basis has to be computed, in order to deal with equations such as \eqref{eqn:TorsionFree}, where both a basis of one-forms and the dual basis of vector fields appear; the member \cpp{dual} of \cpp{Basis} takes care of this.

\smallskip
\cpp{Basis} is based on the standard container \cpp{vector}, part of the Standard Template Library. So, a \cpp{Basis} object can be thought of as a sequence $x_1,\dotsc, x_n$ of \cpp{ex}, where the $x_k$ are linear combinations of simple elements of a given type \cpp{T}, which is a template parameter of \cpp{Basis}. However, \cpp{Basis} differs from \cpp{vector} in the following respects:
\begin{itemize}
\item When elements are added to a \cpp{Basis} object, resulting in a sequence of possibly dependent generators $x_1,\dotsc, x_n$, an elimination scheme is used to remove redunant (dependent) elements. This is done in such a way that the flag $V_1\subseteq \dotsc\subseteq V_n$,
\[V_k=\Span\{x_1,\dotsc, x_k\} \]
is preserved.
\item The first time the member functions \cpp{Components} or \cpp{dual} are called, \cpp{Basis} sets up some internal data according to the following procedure. First, the set
\[S=S(x_1,\dotsc, x_m)=\{\alpha \text{ simple}\mid \langle \alpha,x_k\rangle\neq0 \text{ for some $1\leq k\leq m$}\}\]
is computed; then, the sequence $\{x_1,\dotsc, x_m\}$ that represents the basis internally is enlarged with elements $x_k\in S$, $m< k\leq n$ so that
\[\{x_1,\dotsc, x_n\} \text{ is a basis of }\Span{S}.\]
Then, for each $\alpha\in S$ the components $b_{\alpha j}$ satisfying
\[\alpha=\sum_{j=1}^n b_{\alpha j} x_j\]
are computed; as a matrix, $(b_{\alpha j})$ is the inverse of the matrix with entries $a_{j\alpha}=\langle x_j,\alpha\rangle$.
The dual basis $x^1,\dotsc, x^n$ is then given by
\[x^k=\sum_{\alpha\in S} b_{\alpha k}\alpha.\]
The actual code also take advantage of the fact that with a suitable ordering of the columns, one obtains the block form
\[(a_{j\alpha})=\begin{pmatrix}*&*\\0&I_{n-k}\end{pmatrix}.\]
The dual basis and the $b_{\alpha j}$ are cached internally, until the basis is modified.
\item The member function \cpp{Components} uses the pairing $\langle,\rangle$ to determine the components with respect to $S$, then multiplies by  the matrix $(b_{\alpha j})$ to obtain the components with respect to the basis $x_1,\dotsc x_m$.
\end{itemize}

\smallskip
From the point of view of performance, the operations of \cpp{Basis} have the following complexity, measured in term of the number of multiplications:
\begin{itemize}
\item Given a list $x_1,\dotsc,x_n$ of elements which are known to be linearly independent, one can construct a \cpp{Basis} object with no overhead over \cpp{vector} and no algebraic operation involved.
\item Inserting elements $y_1,\dotsc, y_k$ in a \cpp{Basis} object consisting of elements $x_1,\dotsc, x_n$ has a complexity of $O(m(n+k)^2)$, where \[m=\abs{S(x_1,\dotsc,x_n,y_1,\dotsc,y_k)}.\]
\item Setting up a \cpp{Basis} object consisting of elements $x_1,\dotsc, x_n$ in order that the dual basis or components of a vector may be computed has a complexity of $O((n+m)n^2)$, where \[m=\abs{S(x_1,\dotsc,x_n)}.\]
\item Once the \cpp{Basis} object is so set up, computing the components of a vector has complexity $O(nm)$.
\end{itemize}

Since we work in the symbolic setting, the number of multiplications only gives a rough estimate of execution time. In any case, the above figures
show the importance of being able to define a class that only performs certain operations when it
is needed. Again, this is something which can be achieved naturally in {\CC}, by storing a flag in
each object \cpp{Basis}. Since a \cpp{Basis} object is only modified by invoking certain member
functions, the flag is cleared automatically when the object is modified, effectively invalidating
the cached values of $b_{\alpha j}$.

\begin{remark}
The frames associated to a \cpp{Manifold}, \cpp{RiemannianStructure} or \cpp{Connection} object
are represented by \cpp{Basis} objects; thanks to the mechanisms described in this section, these
frames need not consist of simple elements. In fact, this was the main motivation for introducing \cpp{Basis}.
\end{remark}
\section{An application: Cartan-K\"ahler theory}
In this section we consider, as an application, the problem of determining whether a linear exterior differential system (EDS) is involutive. This problem eventually boils down to computing the rank of certain matrices, but actually writing down these matrices requires a system that supports differential forms and interior products, as well as linear algebra; this makes Wedge particularly appropriate to the task.

Suppose we have a real analytic EDS with independence condition
$(\mathcal{I},\theta_1,\dotsc,\theta_k)$ on a manifold $X$. As a special case, we shall take $X$ to be the bundle of
frames over a $7$-manifold $M$, the $\theta_i$ the components of the tautological form and $\mathcal{I}$  generated by the exterior derivative of the $\Gtwo$-forms $\phi$ and $*\phi$. We wish to determine whether the system is involutive, i.e. every point of $X$ is contained in some submanifold
$S\subset X$ such that all the forms of $\mathcal{I}$ vanish on $S$ and the forms
$\theta_1|_S,\dotsc,\theta_k|_S$ are a basis of the cotangent bundle $T^*S$; we then say that $S$ is an integral manifold for $\mathcal{I}$. In the special case, the EDS turns out to be involutive, proving the existence of local metrics with holonomy $\Gtwo$ (see
\cite{Bryant}).

The general procedure is the following. Complete $\theta_1,\dotsc, \theta_k$ to a basis
$\theta_1,\dotsc,\theta_k,\omega_1,\dotsc,\omega_n$ of $1$-forms on $X$. Suppose that $E\subset
T_xX$ is a $k$-dimensional subspace such that the forms $\theta_1|,\dotsc,\theta_k$ are
independent on $E$; then one can write
\[\omega_i|_E=p_{ij}\theta_j|_E.\]
We say that $E$ is an integral element if the forms of $\mathcal{I}$ restrict to $0$ on $E$. The
$(x,p_{ij})$ are coordinates on the Grassmannian of $k$-planes over $\mathcal{I}$, and by
definition the tangent space of an integral manifold at each point will be an integral element $E$
of the form above. The  integral elements of dimension $n$ form a subset $V_n$ of the
Grassmannian.

Suppose that there exists an integral element $E\subset T_xX$, and $V_n$ is smooth about $E$. By the
Cartan-K\"ahler theory, the (local) existence of an integral  manifold $S$ with $T_xS=E$ is
guaranteed if one can find a flag
\[E_0\subset\dotsc\subset E_k=E, \quad \dim E_k=k,\]
such that
\[\dim V_n= c_0+\dotsc+c_k\;,\]
where $c_j$ is the dimension of the space of polar equations
\[H(E_j)=\{\alpha(e_{1},\dotsc, e_{j},\cdot)\mid\alpha\in\mathcal{I}, e_{l}\in E_j\}\subset T^*_xX.\]
Since we are assuming that $E_j$ is contained in an integral element $E$, $H(E_j)$ has the same dimension as its image
$H^*(E_j)$ under the quotient map
\begin{equation}
 \label{eqn:projection}
T^*_xX\to \frac{T^*_xX}{\langle \theta_1,\dotsc,\theta_k\rangle}.
\end{equation}
We shall say that a system is linear if
\[\mathcal{I}\subset \Lambda^*\Span\{e^1,\dotsc,e^k\}\wedge \Span\{\omega_1,\dotsc, \omega_n\},\]
where the span is over $C^\infty(X)$.
This condition implies that the equations defining $V_n$ are affine in the $p_{ij}$. In addition, one has $H^*(E_j)=H^*(E_j')$, where $E_j'$ is the horizontal projection of $E_j$ determined by the splitting
\[\langle \theta_1,\dotsc,\theta_k,\rangle \oplus\langle\omega_1,\dotsc,\omega_n\rangle.\]

One can compute the dimension of $V_n$ and the $c_i$ with the surprisingly short program of Fig. \ref{fig:Cartan}. This program represents the frame bundle $P$ as a parallelizable manifold with a frame $(\theta_i,\omega_{jk})$, where the $\theta_i$ are the components of the tautological form, and the $\omega_{jk}$ are the components of a torsion-free connection form, related by the structure equation
\[de^i=e^j\wedge\omega_{ij}.\]
The member function \cpp{GetEquationsForVn} simply replaces every occurrence of $\omega_{jk}$ with $p_{jk}$ in terms of the coordinates $p_{ij}$ of the Grassmannian, and imposes that the coefficients of the resulting forms are zero. The member function \cpp{GetReducedPolarEquations}, on the other hand, calls itself recursively to obtain all the forms $e_{i_1}\hook\dotsc\hook e_{i_l}\hook \alpha$, where $\alpha$ ranges in a basis of $\mathcal{I}$ and $1\leq i_1<\dotsc<i_l\leq j$. When the resulting degree is one, the function stops the recursion and stores the differential one-form, applying a substitution corresponding to the projection~\eqref{eqn:projection}.
\begin{figure}[tb]
\smallskip
\begin{verbatimtab}
struct EDS :  public ConcreteManifold, public Has_dTable  {
    lst moduloIC;
    int dim;

    EDS(int dimension) : ConcreteManifold(dimension*(dimension+1)) {
        dim=dimension;
        for (int i=1;i<=dim;++i) {
            ex x;
            for (int j=1; j<=dim;++j)
                x+=e(j)*e(i*dim+j);
            Declare_d(e(i),x);
            moduloIC.append(e(i)==0);
        }
    }

    template<typename Container>
    void GetEquationsForVn(Container& container, lst I) {
        lst substitutions;
        for (unsigned i=dim+1; i<=dim*(dim+1);++i) {
            ex ei;
            for (unsigned j=1;j<=dim;++j)
                ei+=symbol("p"+ToString(i,j))*e(j);
            substitutions.append(e(i)==ei);
        }
        for (lst::const_iterator i=I.begin(); i!=I.end();++i)
            GetCoefficients(container,i->subs(substitutions));
    }

    template<typename Container>
    void GetReducedPolarEquations(Container& container, ex form, int j) {
        if (degree<DifferentialForm>(form)==1)
            container.push_back(form.subs(moduloIC));
        else if (j>0) {
            GetReducedPolarEquations(container, form,j-1);
            GetReducedPolarEquations(container, Hook(e(j),form),j-1);
        }
    }
};
\end{verbatimtab}
\caption{Using Wedge to determine if a linear EDS on the frame bundle is involutive}
\label{fig:Cartan}
\end{figure}

The code of Fig.~\ref{fig:Cartan} is for a general linear EDS on the bundle of frames. In order to obtain a result for the special case of $\Gtwo$, we can use the following:

\begin{verbatimtab}
struct G2 : EDS {
    G2() : EDS (7) {
        lst I;
        I=d(ParseDifferentialForm(e(),"567-512-534-613-642-714-723")),
          d(ParseDifferentialForm(e(),"1234-6712-6734-7513-7542-5614-5623"));
        Basis<symbol> A;
        GetEquationsForVn(A,I);
        for (int j=0;j<7;++j)
        {
            Basis<DifferentialForm> V;
            for (lst::const_iterator i=I.begin();i!=I.end();++i)
                GetReducedPolarEquations(V,*i,j);
            cout<<V.size()<<endl;
        }
        cout<<A.size()<<endl;
    }
};
\end{verbatimtab}

The output shows that the $c_i$ are
\[0,0,0,1,5,15,28\]
whereas the codimension of $V_n$ is $49$; thus, Cartan's test applies and the system is involutive.

Notice that the equations on the Grassmannian are stored in a \cpp{Basis<symbol>} object, and so the codimension of $V_7$ is just the size of the basis. This makes sense because we know that the equations defining $V_7$ are linear in the $p_{ij}$; however, the container \cpp{A} is passed as a template parameter to \cpp{GetEquationsForVn}. This means that if we were dealing with an EDS  where the equation are affine in the $p_{ij}$ rather than linear, e.g. the EDS associated to nearly-K\"ahler structures in six dimensions, it would suffice to modify the code by declaring \cpp{A} to be of type \cpp{AffineBasis<symbol>}. The general philosophy, used extensively in Wedge, is that of writing ``generic'' template functions, without making specific assumptions on the nature of the arguments, so that the type of the arguments determines the precise behaviour of the code.

\begin{remark}
Concerning the non-linear case, Wedge also provides a class template \cpp{PolyBasis}, analogous to
\cpp{Basis}, where polynomial equations can be stored, relying on {\hbox{\rm C\kern-.13em
o\kern-.07em C\kern-.13em o\kern-.15em A}} for the necessary Gr\"obner basis computations (see
\cite{CocoaSystem}). Notice however that the simple program of Fig.~\ref{fig:Cartan} assumes
linearity when computing the polar equations.
\end{remark}

\bibliographystyle{plain}
\bibliography{wedge}

\end{document}